\documentclass[a4paper,12pt]{amsart}
\usepackage{a4,amsbsy,amscd,amssymb,amstext,amsmath,latexsym,oldgerm,enumerate}

\newtheorem{thmgl} {Theorem}              
\newtheorem{propgl}{Proposition}
\newtheorem{lemgl} {Lemma}
\newtheorem{cornn}{Corollary}

\newcommand{\mf}{\mathfrak}
\newcommand{\mc}{\mathcal}
\newcommand{\mb}{\mathbb}

\newcommand{\nts}{\negthinspace}     

\newcommand{\ncd}{\nts\cdot\nts}

\newcommand{\sm}{\setminus}         

\newcommand{\ot}{\otimes}           

\newcommand{\Der}{{\rm Der}}

\newcommand{\Ker}{{\rm Ker}}

\newcommand{\codim}{{\rm codim}}
\newcommand{\sgn}{{\rm sgn}}
\newcommand{\tr}{{\rm tr}}

\newcommand{\gr}{{\rm gr}}   

\newcommand{\g}{{\mf{g}}}
\newcommand{\h}{{\mf{h}}}
\renewcommand{\t}{\mf{t}}
\newcommand{\n}{\mf{n}}

\newcommand{\z}{\mf{z}}
\newcommand{\gl}{\mf{gl}}
\newcommand{\spl}{\mf{sl}}
\newcommand{\psl}{\mf{psl}}
\newcommand{\pgl}{\mf{pgl}}

\newcommand{\Ad}{{\rm Ad}}              
\newcommand{\Lie}{\rm{Lie}}
\newcommand{\GL}{{\rm GL}}
\newcommand{\SL}{{\rm SL}}

\begin{document}
\title{Zassenhaus varieties of  general linear Lie algebras}
\author{Alexander Premet and Rudolf Tange\\
\medskip}
\begin{abstract}
Let $\g$ be a Lie algebra over an algebraically closed field of characteristic
$p>0$
and let $U(\g)$  be the universal enveloping   
algebra of $\g$.
We prove in this paper that
for $\g=\gl_n$ and $\g=\spl_n$  the centre of
$U(\g)$ is a unique factorisation domain and 
its field of fractions  
is rational. For $\g=\spl_n$ our argument requires  
the assumption that $p\nmid n$ while for $\g=\gl_n$ it works for any $p$.  
It turned out that our two main results are
closely related to each other. 
The first one confirms in type ${\rm A}$ a recent conjecture of
 A.~Braun and C.~Hajarnavis while the second  answers a question of J.~Alev.
\end{abstract}
\maketitle
\section{Introduction}
Let $K$ be an algebraically closed field of characteristic $p>0$. 
In this note $G$  denotes a connected reductive 
$K$-group  with Lie algebra $\g$. Mostly we will be in the situation where
$G=\GL_n(K)$ or $G=\SL_n(K)$ and $p\nmid n$. Let $x\mapsto x^{[p]}$ denote the canonical
$p$th power map on $\g$ equivariant under the adjoint action of $G$.

The group $G$ acts on $U$ as algebra 
automorphisms. This action extends the adjoint action of $G$ on $\g$, hence preserves 
the standard filtration $(U_i)_{i\ge 0}$ of $U$. 
The associated graded algebra $\gr(U)=S(\g)$ is a domain, so $U$ has no zero 
divisors. The centre $Z$ of $U$ is therefore a filtered $K$-algebra, a domain, and a 
filtered $G$-module.

Let ${\mc Q}={\mc Q}(\g)$ be the field of fractions of $Z$.
By a classical result of Zassenhaus, 
$Z$ is Noetherian and integrally closed in 
${\mc Q}$; see \cite{Z}. Moreover, 
$\mathrm{tr.\,deg}_K{\mc Q}=\dim \g$ and
the localisation
$\mc{D}(\g):={\mc Q}\ot_{Z}U(\g)$ is a central division 
algebra over $\mc Q$ of dimension $N^2$ where $N$ 
is the maximal dimension of irreducible $\g$-modules. When
$G=\GL_n(K)$ or $G=\SL_n(K)$ we have $N=p^{n(n-1)/2}$; see
\cite{Jan2} for example. The  maximal spectrum 
${\mc Z}$ of the algebra $Z$ 
is called the {\it Zassenhaus variety} of $\g$. 
By the above discussion,
the variety $\mc Z$ is affine, irreducible and normal. Furthermore,
$\dim\,{\mc Z}=\dim\,\g$.  It is proved in \cite{BG} that
under rather mild assumptions on $p$ the singular points of $\mc Z$
are exactly the maximal ideals $\mf m$ for which
$(Z/{\mf m}Z)\ot_{Z}U$ is not isomorphic to the matrix algebra $\text{Mat}_N(K)$.

At present very little is known about the division algebra ${\mc D}(\g)$ and its
class in 
the Brauer group of ${\mc Q}$. 
In order to get started here it will be important to
address the following question posed to the first author by Jacques Alev:

\smallskip
\noindent
{\bf Question} (J.~Alev). {\it Is it true that $\mc Q$ is $K$-isomorphic to the field 
of rational functions  $K(X_1,\ldots, X_m)$ with $m=\dim \g$? In other words, is it true 
that the Zassenhaus variety $\mc Z$ is rational?}

 \smallskip 
Until now the answer to this question was known  only in the simplest
case $\g=\spl_2$. Another interesting question
related to $\mc Z$ was recently raised  in 
\cite{BHa}  and answered positively for  
$\g=\spl_2$ (mild characteristic restrictions may apply): 

\smallskip
\noindent
{\bf Conjecture} (A.~Braun and C.~Hajarnavis). {\it The centre of $U(\g)$ is a 
unique factorisation domain.}

\smallskip

Similar problems can be raised in the characteristic zero case as well.
Here one has to replace $U(\g)$ by the quantised enveloping algebra 
$U_{\epsilon}(\g_{\mb C})$
(without divided powers) at a root of unity $\epsilon\in\mb C$; 
see \cite{BHa} for more detail.

The main result of this paper is the following theorem which solves
both problems in the modular case for $\g=\gl_n$ and for $\g=\spl_n$ with
$p\nmid n$:

\smallskip

\noindent
{\bf Theorem.} {\it If $\g=\gl_n$ or $\g=\spl_n$ and
$p\!\nmid \!n$, then the centre of $U(\g)$ is a unique factorisation 
domain and its field of fractions is rational.}

\smallskip

One expects this result to extend to the
Lie algebras $\g$ isomorphic to 
$\spl_n$, $\pgl_n$ and $\psl_n$ with $p\!\mid\! n$. However,
to obtain such an extension by our methods
one would  need an explicit 
description of the invariant 
algebra $S(\g)^\g$, which is currently unavailable.
As for the Lie algebras 
of other types, both problems remain open and 
new ideas are required here.

\smallskip

\noindent
{\bf Acknowledgement.} We would like to thank Serge Skryabin 
for drawing our attention to Proposition~1. It allowed us to simplify
our original proof of  the main theorem.

\section{Preliminaries}

\noindent{\bf 2.1.}\quad
Given an element $x$ of a commutative
ring $S$ we denote by $(x)$ the ideal of $S$ generated by $x$. Recall that $x$ is called
{\it prime} if $(x)$ is a prime ideal of $S$.

 Let $A$ be an associative ring with an
ascending filtration $(A_i)_{i\in\mb{Z}}$. If $I$ is a two sided ideal of $A$, then the abelian group $I$ and the
ring $A/I$ inherit an ascending filtration from A and we have an embedding
$\gr(I)\hookrightarrow\gr(A)$ of graded abelian groups. If we identify $\gr(I)$ with a
graded subgroup of the graded additive group $\gr(A)$ by means of this embedding, then
$\gr(I)$ is a two sided ideal of $\gr(A)$ and there is an isomorphism
$\gr(A/I)\cong \gr(A)/\gr(I)$; see \cite{Bou2}, Chapter~3,  \S~ 2.4.

Now assume that
$\bigcup_iA_i=A$ and $\bigcap_iA_i=\{0\}$. For a nonzero $x\in A$ we define
$\deg(x):=\min\{i\in\mb{Z}\,|\,\,x\in A_i\}$ and $\gr(x):=x+A_{k-1}\in\gr(A)^k=A_k/A_{k-1}$ where
$k=\deg(x)$. If $\gr(A)$ has no zero divisors, then the same holds for $A$ and we have for
$x,y\in A\sm\{0\}$ that
$\deg(xy)=\deg(x)+\deg(y)$, $\gr(xy)=\gr(x)\gr(y),$ and $\gr((x))=(\gr(x))$. We mention
for completeness
that if $A=\bigoplus_{n\in\mb{Z}} A^n$ is a graded ring, then $(A_n)_{n\in\mb{Z}}=
\big(\sum_{k\le n}A^k\big)_{n\in\mb{Z}}$ defines an ascending filtration of $A$ with the
two properties
mentioned above and $A\cong\gr(A)$ as algebras.
\smallskip

\noindent
{\bf 2.2.}\quad
The {\it $p$-centre}
$Z_p$ of $U$ is defined as the subalgebra of $U$ generated by all elements
$x^p-x^{[p]}$ with $x\in \g$. It is  well-known
(and easily seen) that $Z_p\subseteq Z$ is a polynomial algebra in $x_i^p-x_i^{[p]}$
where $\{x_i\}$ is any basis of $\g$. For a vector space $V$ over
$K$ the {\it Frobenius twist}
$V^{(1)}$ of $V$ is defined as the vector space over $K$ with the
same additive group as $V$ and with scalar multiplication given by
$\lambda\ncd x=\lambda^{1/p}\,x$.
Note
that the linear functionals and the polynomial functions on $V^{(1)}$ are the $p$-th powers
of
those of $V$. The Frobenius twist of a $K$-algebra is defined similarly (only the scalar
multiplication is modified). Following \cite{Kry} we define $\eta\,\colon\,S(\g)^{(1)}\to Z_p$ by
setting
$\eta(x)=x^p-x^{[p]}$ for all $x\in\g$; see also \cite{MiRu}. This is a $G$-equivariant algebra
isomorphism, hence it restricts to an algebra isomorphism
$$\eta\,\colon\,\big(S(\g)^G\big){^{(1)}}=\big(S(\g)^{(1)}\big)^G
\stackrel{\sim}{\rightarrow} \,Z_p^G.$$ We have
$\gr(\eta(x))=x^p$ for all $x\in\g\sm\{0\}$. Furthermore
the associated graded algebra of the filtered algebra $Z_p\subset U$
is $G$-equivariantly isomorphic to
the graded subalgebra $S(\g)^p$ of $S(\g)$.

\smallskip
\noindent {\bf 2.3.}\quad In the remainder of this note we assume
that $G=\GL_n(K)$ or $G=\SL_n(K)$ and $p\nmid n$. In this case
Theorem~1.4 in \cite{FrPa} shows that the filtered $G$-modules
$U(\g)$ and $S(\g)$ are isomorphic (the isomorphism in \cite{FrPa}
is obtained by composing the Mil'ner map $\phi\,\colon\, U\to \,
S(U)$ with a $G$-equivariant projection from $U$ onto $\g$).
Consequently, each $G$-module $U_{n-1}$ has a $G$-invariant direct
complement in $U_{n}$. This implies that the associated graded
algebras of $U^G$ and $Z$ are isomorphic to $S(\g)^G$ and
$S(\g)^\g$, respectively.

The trace form $\beta\,\colon\,\gl_n\times\gl_n\to K$
associated with the vector representation of $\GL_n(K)$ is  nondegenerate
and the same holds for its restriction
to $\spl_n$ as $p\nmid n$. Let $\theta:S(\g^*)\to
S(\g)$ denote the
$G$-equivariant algebra isomorphism induced by $\beta$ (it takes
$f\in\g^*$ to a unique $x\in\g$ such that $f(y)=\beta(x,y)$ for all $y\in\g$).

Let $\h$ be the subalgebra of all diagonal matrices in $\gl_n$ and $\h'=\h\cap\spl_n$.
Let
$\n^+$ (respectively, $\n^-$) be the subalgebra of all strictly upper (respectively,
lower) triangular matrices in $\g$. To unify notation we set
$\t=\h$ if $\g=\gl_n$ and $\t=\h'$ if $\g=\spl_n$.
Then  we have $\g=\n^{-}\oplus \t \oplus\n^+$.
Also, $\t=\Lie\,T$ where $T$ is the group of all diagonal matrices in $G$.
Furthermore, $\t$ is the orthogonal
complement
to $\n^-\oplus\n^+$  with respect to $\beta$. The Weyl group action  induced by
the adjoint action of the normaliser $N_G(T)$ on $\t$ is nothing but the restriction to
$\t$ of the permutation action of the symmetric group ${\mf S}_n$ on the space of
diagonal matrices $\h$.

In \cite{KW}, Theorem 4, Kac and Weisfeiler proved
that a modular version of the
Chevalley restriction theorem holds for the coadjoint action of any
simple,  simply connected algebraic $K$-group.
Their arguments are known to work for  all connected reductive $K$-groups
with simply connected derived subgroups; see \cite{Jan2}, Sect.~ 9
for example. In particular, they apply to our group $G$.
Since $\theta\,\colon\,K[\g]\to K[\g^*]$ is a $G$-equivariant algebra isomorphism,
Theorem~4 in \cite{KW} implies that
the restriction map $K[\g]\to K[\t]$ induces
an algebra isomorphism $K[\g]^{G}\stackrel{\sim}{\longrightarrow} \,K[\t]^{{\mf S}_n}$.

\smallskip
\noindent
{\bf 2.4.}\quad
For $1\le i\le n$ define $s_i\in K[\gl_n]^{\GL_n}$ by setting
$s_i(x)=\tr\,(\wedge^i\, x)$ for all
$ x\in\gl_n,$ where $\wedge^i\, x$ is the $i$th exterior power of $x$.
Then
$$\chi_x(X)=X^n+\sum_{i=1}^n (-1)^i s_i(x)X^{n-i}$$ is the characteristic
polynomial of $x$. Let $\{e_{i,j}\,|\,1\le i,j\le n\}$ be the basis of $\gl_n$ consisting of
the matrix units and let
$\{\xi_{ij}\,|\,1\leq i,j\leq n\}$ be the corresponding dual basis in $\gl_n^*$. To ease notation
identify each $\xi_{ii}$ with its restriction to the diagonal subalgebra $\h$.
For $1\le i\le n$ the restriction of
$s_i$ to $\h$ is then the $i$th elementary symmetric function
$\sigma_i$ in $\xi_{11},\,\xi_{22},\,\ldots, \xi_{nn}$.
By
the theorem on symmetric functions, $\sigma_1,\ldots,\sigma_n$
are algebraically independent and
generate the invariant algebra $K[\h]^{{\mf S}_n}$. Our discussion in 2.3
now shows that the $s_i$'s  are
algebraically independent and
generate the invariant algebra $K[\gl_n]^{\GL_n}$.

Suppose $p\!\nmid\! n$. Given a polynomial function $f$ on $\gl_n$ we
denote by $f'$ its restriction to $\spl_n$. The span of all
$\xi_{ii}-\xi_{jj}$ is an ${\mf S}_n$-invariant direct complement to the line $K \sigma_1$ in
$\h^*$, hence the $K$-subalgebra  generated by all $\xi_{ii}-\xi_{jj}$ is an
${\mf S}_n$-invariant direct
complement to the ideal of $K[\h]$ generated by $\sigma_1$. From this it is immediate
that the restriction map $K[\h]\to K[\h']$ induces an epimorphism
$K[\h]^{{\mf S}_n}\twoheadrightarrow K[\h']^{{\mf S}_n}$ whose kernel is the ideal of
$K[\h]^{{\mf S}_n}$ generated by $\sigma_1$.
Since the subalgebra of $K[\h]^{{\mf S}_n}$ generated by
$\sigma_2,\ldots,\sigma_n$ is a direct complement in $K[\h]^{{\mf S}_n}$ to this ideal, we deduce that
the restrictions $s_2'\vert_{\h'},\ldots,s_n'\vert_{\h'}$
are algebraically independent and generate $K[\h']^{{\mf S}_n}$.
But then $s_2',\ldots,s_n'$ are
algebraically independent and generate the invariant algebra $K[\spl_n]^{\SL_n}$ by
our discussion in 2.3.

Under the
$G$-equivariant isomorphism $\theta\,\colon\,
S(\g^*)\stackrel{\sim}{\longrightarrow} \,
S(\g)$ and the induced ${\mf S}_n$-equivariant
isomorphism $S(\t^*)\stackrel{\sim}{\longrightarrow} \,
 S(\t)$, the restriction map
$S(\g^*)\to S(\t^*)$ corresponds to
the projection homomorphism $\Phi:S(\g)\to S(\t)$
defined as follows: if we identify $S(\g)$ with
$S(\n^-)\ot S(\t)\ot S(\n^+)$, then $\Phi(x\ot h\ot y)=x^0hy^0$ where $f^0$
denotes the zero degree part of $f\in S(\g)$. By the above, $\Phi$ induces an
algebra isomorphism $S(\g)^G\cong S(\t)^{{\mf S}_n}$.

\smallskip
\noindent {\bf 2.5.}\quad
In \cite{KW}, Kac and Weisfeiler also proved  a
noncommutative  version of the Chevalley restriction theorem. Again the
arguments in \cite{KW} are known to generalise to
all connected reductive $K$-groups
with simply connected derived subgroups.  In particular, they apply to our group $G$.

Let $\Psi\,\,\colon
U=U(\n^-)\ot U(\t)\ot U(\n^+)\longrightarrow\, U(\t)=S(\t)$ be the linear map taking
$x\ot h\ot y$ to $x^0hy^0$, where $u^0$ denotes the scalar part of $u\in U$ with respect
to the decomposition $U=K1\oplus\, U_+$ where $U_+$ is the
augmentation ideal of $U$.  The restriction of $\Psi$
to $U^{N_G(T)}$ is an algebra homomorphism.
Define the shift homomorphism $\gamma\,\colon\,S(\t)\to S(\t)$
by setting
$\gamma(h)=h-\rho(h)$ for all $h\in\t$, where $\rho$ denotes the half sum of
the differentials
of the positive roots. It follows from \cite{KW}, Theorem~1, that
$\gamma\circ\Psi$ induces an algebra isomorphism between
$U^G$ and  $S(\h)^{{\mf S}_n}$. As a consequence,
$U^G$ is a polynomial algebra in $\dim\,\t$ variables.

Using the descriptions of $\Phi$ and $\Psi$ and a PBW-basis it follows that for  $x\in U\sm\{0\}$ with
$\Phi(\gr(x))\neq 0$ we have $\Psi(x)\neq 0$ and $$\gr(\gamma(\Psi(x)))\,=\,\gr(\Psi(x))\,=\,
\Phi(\gr(x)).$$
By the injectivity of the restriction of
$\Phi$ to $S(\g)^G$, the displayed  equalities hold for all
$x\in U^G$. Thus we can deduce the injectivity of $\gamma\circ\Psi\,\colon\,
U^G\longrightarrow \,
S(\t)^{{\mf S}_n}$
from that of $\Phi\,\colon\,S(\g)^G
\longrightarrow \,S(\t)^{{\mf S}_n}$. The same applies to the surjectivity; see the
proof of Proposition~2.1 in \cite{Ve}.
\smallskip
\section{Invariants for the Lie algebra}\label{g-invariants}
\noindent
{\bf 3.1.}\quad
The aim of this section is to put together all results on Lie algebra invariants
that will be in use later on. These results are
mostly known but
their proofs  are spread over the
literature (and folklore); see
\cite{Ve}, \cite{KW}, \cite{FrPa}, \cite{D}, \cite{BrGo}, \cite{Skry},
\cite{Jan3}, Sect.~7, and the references therein.

Given $x\in\g$ we denote by $\z_{\g}(x)$ the centraliser of $x$ in $\g$.
An element $x\in\g$ is called {\it regular}
if $\dim\,\z_\g(x)=\dim \t$.
It is well-known and not hard to see that $\dim\,\z_{\g}(x)\ge \dim \t$
for all $x\in\g$.\footnote{As in the group case,
take a Borel subgroup $B$ of $G$ with
$x\in{\Lie}(B)$ and consider the morphism $B\to{\Lie}(B,B)$
sending $g\in B$ to  $({\Ad}\,g)(x)-x \in {\Lie}(B,B)$; see  \cite{St}, page~1.}
Moreover,
the set $\g_{\text{reg}}$ of all regular elements
in $\g$ is  nonempty and Zariski open in $\g$.
Furthermore,  Linear Algebra
shows that $x$ is regular in $\gl_n$ if and only if
the minimal polynomial of $x$ equals $\chi_x(X)$,
which
happens if and only if the column space $K^n$
is a cyclic $K[x]$-module.

The first result we need is a modular version of Kostant's differential criterion of
regularity. It is is essentially due to Veldkamp \cite{Ve}.

\begin{lemgl}
For $x\in\gl_n$ the following  are equivalent:
\begin{enumerate}
\item[(1)] the element $x$ is regular;
\item[(2)]  the differentials $d_xs_1,\ldots,d_xs_n$ are linearly independent.
\end{enumerate}
\end{lemgl}

\begin{proof}
That the independence of $d_xs_1,\ldots,d_xs_n$ implies
the regularity of $x$ is proved in  \cite{Ve}, Sect.~7.
The proof requires a lemma on the invariant algebra
$K[\g]^G$ (Lemma 7.2 in \cite{Ve}), the fact that the semisimple irregular elements
of $\g$ form
a dense subset in $\g\setminus\g_{\text{reg}}$
(Proposition 4.9 in \cite{Ve}), and a result from
\cite{Bou1} (Prop.~6, Chap.~5, Sect.~ 5.5). All these are valid for $\g=\gl_n$.

That  the regularity of $x$ implies the
independence of $d_xs_1,\ldots,d_xs_n$ is much easier to prove.
Given
${\bf a}=(a_1,\ldots,a_n)\in K^n$ we set
$$x_{\bf a}\,=\,
\left(\begin{array}{ccccc}
a_1&a_2&\cdots&a_{n-1}&a_n\\
1&0&\cdots&0&0\\
0&1&\cdots&0&0\\
\vdots&\vdots&\ddots&\vdots&\vdots \\
0&0&\cdots&1&0\end{array}\right).$$
Each  $x_{\bf a}$ is regular  in
$\gl_n$ as the minimal polynomial of
$x_{\bf a}$ equals
$X^n-\sum_{i=1}^n a_iX^{n-i}$. The set  ${\mc S}=\{x_{\bf a}\,|\, {\bf a}\in K^n\}$
is an $n$-dimensional affine subspace in $\gl_n$ through the point
$x_{\bf 0}$. The restriction to
$\mc S$ of the morphism $x\longmapsto(s_1(x),\ldots,s_n(x))$ is an isomorphism
of $\mc S$ onto ${\mb A}^n$. From this it is immediate that
the differentials
$d_xs_1,\ldots,d_xs_n$ are linearly independent for all $x\in\mc S$.
On the other hand, every matrix $x$ whose minimal polynomial equals $\chi_x(X)$
is similar to a matrix from $\mc S$. Hence
these differentials are independent for all  regular $x$.
\end{proof}
\noindent {\bf 3.2.}\quad Now we look at the regular elements in
$\spl_n$. Recall the notational conventions of 2.4. It is
immediate from the definition that  $x\in\gl_n$ is regular if and
only so is $x+\lambda I_n$  for any $\lambda\in K$.
\begin{cornn}
Suppose $p\nmid n$. For $x\in\spl_n$
the following are equivalent:
\begin{enumerate}[(1)]
\item the element $x$ is regular in $\spl_n$;
\item the element $x$ is regular in $\gl_n$;
\item the differentials $d_xs'_2,\ldots,d_xs'_n$ are linearly independent.
\end{enumerate}
\end{cornn}
\begin{proof}
We have $\z_{\gl_n}(x)\,=\,\z_{\spl_n}(x)\oplus KI_n$.
This shows that (1) and (2) are equivalent. The differentials
$d_xs_1,\ldots,d_xs_n$ are independent if and only so are the restrictions of
$d_xs_2,\ldots,d_xs_n$ to $\spl_n$, the kernel of $d_xs_1=s_1$.
The equivalence of (2) and (3) now follows from Lemma~1.
\end{proof}
\noindent
{\bf 3.3.}\quad As mentioned in the Introduction, our proof of the main theorem
will rely on the following proposition communicated to us by S.~Skryabin.
We were unable to trace this result in the literature. Although
it resembles strongly one of the basic facts of the invariant theory  of groups,
it also captures some essential features of the invariant theory
of restricted Lie algebras.

Recall that the coordinate algebra $K[V]$ of
a finite dimensional vector space $V$ over $K$
is a unique factorisation domain. The algebra
$K[V]\,\cong \,\bigoplus_{i\ge 0}S^i(V^*)$ is graded
and $\gl(V)$ acts on $K[V]$ as
homogeneous derivations of degree $0$. Therefore
$K[V]^p\subseteq K[V]^{\mc L}$.
\begin{propgl}
Let $L$ be a Lie algebra such that $L=[L,L]$ and let
$V$ be a finite dimensional $L$-module. Then the invariant algebra
$K[V]^L$ is a unique factorisation domain and the irreducible
elements
of $K[V]^L$ are the $p$th powers of the irreducible elements of
$K[V]$ not invariant under $L$
and the irreducible elements of $K[V]$ contained in  $K[V]^L$.
\end{propgl}
\begin{proof}
Let $f$ be a nonzero element in $K[V]^L$ and suppose $f=f_1f_2$ where
$f_1,f_2\in K[V]$ are coprime of positive degree. Let $x$
be any element in $L$. Since $(x\cdot f_1)f_2=-f_1(x\cdot f_2)$,
the uniqueness of prime factorisation in $K[V]$ implies that $f_2$
divides $x\cdot f_2$. As $\deg(x\cdot f_2)\le \deg f_2$ it must be that $x\cdot f_2=\chi(x)f_2$
for some $\chi(x)\in K$. The map $\chi\,\colon\,L\to K$ is a character of $L$. As
$L=[L,L]$, it must be that $\chi=0$. This shows that $f_1,f_2\in K[V]^L$.
Now suppose $f=g^n$ for some $n\in \mb N$. Write $n=sp+r$  with
$s,r\in {\mb Z}_+$ and $0\le r<p$. Then
$0=x\cdot f=ng^{n-1}(x\cdot g)$. For $r\ne 0$ this yields $g\in K[V]^L$, while for
$r=0$ we have  $f=(g^p)^s$ with $g^p\in K[V]^L$.

This shows that any irreducible
element in $K[V]^L$ is either an irreducible element of $K[V]$ invariant
under $L$
or a $p$th power of an irreducible element in $K[V]\setminus K[V]^L$.
Now the
unique factorisation property of
$K[V]^L$ follows from that of $K[V]$.
\end{proof}
\noindent
{\bf 3.4.}\quad
Let $X$ be an affine algebraic variety defined over $K$,  and  let $\mc L$
be a finite dimensional
restricted Lie algebra together with a restricted homomorphism ${\mc L}\to\Der_K\,
K[X]$. Define ${\mc L}_x$ to be the stabiliser of the maximal ideal
${\mf m}_x$ of $K[X]$ corresponding to a point $x\in X$. Following \cite{Skry},
Sect.~5,  we put $$c_{\mc L}(X):=\,\max_{x\in X}\,
\codim_{\mc L}\,{\mc L}_x.$$
In the situation of 3.3,
where $X=V$ is a finite dimensional restricted $\mc L$-module, it is
easy to see that
${\mc L}_x=\{l\in{\mc L}\,|\,\,l(x)=0\}$ for every $x\in V$.
\begin{lemgl}
We have $K[\gl_n]^{\gl_n}\,=\,K[\gl_n]^{\spl_n}$ for all $n\in\mb N$. Moreover,
$K[\gl_n]^{\gl_n}$ is a unique factorisation domain and the irreducible
elements
of $K[\gl_n]^{\gl_n}$ are the $p$th powers of the irreducible elements of
$K[\gl_n]$ not invariant under $\gl_n$
and the irreducible elements of $K[\gl_n]$ contained in $K[\gl_n]^{\gl_n}$.
\end{lemgl}
\begin{proof} 1. For $p\nmid n$ the first
part of the statement is obvious as
$\gl_n=\spl_n\oplus KI_n$. To tackle it in the general case we
recall our notation in 2.3 and
set $V=\gl_n$.
It follows from our remarks above that
$(\gl_n)_x=\z_{\gl_n}(x)$ for all $x\in V$. So
the discussion in 3.1 yields that $c_{\gl_n}(V)=n^2-n$.
Let $h$ be a regular element of $\gl_n$ contained in $\h$.
Then we have $(\gl_n)_h=\z_{\gl_n}(h)=\h$ and $\gl_n=\spl_n+(\gl_n)_h$.
But then $K[\gl_n]^{\gl_n}\,=\,K[\gl_n]^{\spl_n}$ in view of \cite{Skry}, Corollary~5.3.

\smallskip

\noindent
2. The second part of the statement follows immediately from Proposition~1 if $(p,n)\neq (2,2)$,
since then, as is well-known, $\spl_n$ is perfect.
To establish it in general we will
slightly modify our arguments in the
proof of Proposition~1. If for $f\in K[V]^{\gl_n}$ we have $f=f_1\,f_2$ with
$f_1,f_2\in K[V]$ coprime, then as in that proof $x\cdot f_2=\chi(x)\,f_2$ for all $x\in\gl_n$.
The character $\chi:\gl_n\to K$ must vanish on $[\gl_n,\gl_n]=\spl_n$.
But then $f_1, f_2\in K[V]^{\gl_n}$, by part~1 of this proof. The rest of the proof
of Proposition~1 applies in our present situation, and the result follows.
\end{proof}
\noindent
{\bf 3.5.}\quad
The statement below is known but  we wanted to streamline its proof
by employing the relationship between filtered and graded algebras  in a more
systematic way.
Assertion  (iv)  is often referred to as Veldkamp's
theorem; see \cite{Ve} Theorem 3.1.
\begin{propgl}\label{prop.g-invariants}
Let $m$ be the rank of
$\g$, i.e. $m=\dim \mf t$,
and put $(t_1,\ldots,t_m)=(s_1,\ldots,s_n)$ for $\g=\gl_n$ and
$(t_1,\ldots,t_m)=(s_2',\ldots,s_n')$ for $\g=\spl_n$. Define $u_i\in U^G$
by setting
$u_i=\big((\gamma\circ\Psi)^{-1}\circ\Phi\big)\big(\theta(t_i)\big)$. Then the following
hold: \begin{enumerate}[(i)]

\item The set $\g\setminus\g_{\mathrm{reg}}$ is Zariski closed of pure codimension
$3$ in $\g$.

\item $K[\g]^\g$ is a free $K[\g]^p$-module with basis
$\{t_1^{k_1}\cdots t_m^{k_m}\,|\,\,1\leq k_i<p\}$.

\item $S(\g)^\g$ is a free $S(\g)^p$-module with basis
$\{\theta(t_1)^{k_1}\cdots \theta(t_m)^{k_m}\,|\,\,1\leq k_i<p\}$.

\item $Z$ is a free $Z_p$-module with basis $\{u_1^{k_1}\cdots u_m^{k_m}\,|\, \,
1\leq k_i<p\}$.
\end{enumerate}
\end{propgl}
\begin{proof}
(i) The first assertion
 is proved in \cite{Ve} Theorem 4.12. The arguments there also apply to
$\g=\gl_n$.

\noindent
(ii) By Lemma~1, Corollary~1 and (i), the Zariski closed subset of $\g$
consisting of all $x$ for which the differentials
$d_xt_1,\ldots,d_xt_m$
are linearly dependent has codimension 3 in $\g$.
The second assertion now follows from Theorem 5.4 in \cite{Skry}
applied to the variety $X=\g$. Arguing as in the proof
of Lemma~2 one observes that
$c_\g(X)=n^2-n$ in our case.  Therefore, $\dim X-c_\g(X)=m$.

\noindent
(iii) The third assertion follows immediately from part (ii)
in view of the isomorphism  $\theta\colon\,K[\g]\,\stackrel{\sim}{\longrightarrow}\,
S(\g)$.

\noindent
(iv) Recall from  2.2 and 2.3 that the associated graded algebras of
$Z$, $U^G$ and $Z_p$ are $S(\g)^\g$, $S(\g)^G$ and $S(\g)^p$, respectively.
By our remarks in  2.3 and 2.5
we have $\theta(t_i)=\gr(u_i)$. The fourth
assertion  now follows from part (iii) by a standard induction argument;
see the proof of Theorem~3.1 in \cite{Ve} for more detail.
\end{proof}
{\bf Remarks.} 1.
It follows from Proposition~2
that the bases in (ii), (iii), (iv)  are also
bases of $K[\g]^G$, $S(\g)^G$ and $U^G$ over  $\big(K[\g]^p\big)^G$, $\big(S(\g)^p\big)^G$
and $Z_p^G$, respectively.
This implies that
$K[\g]^\g\cong K[\g]^p\ot_{(K[\g]^p)^G} K[\g]^G$,
$S(\g)^\g\cong S(\g)^p\ot_{(S(\g)^p)^G} S(\g)^G$
and $Z\cong Z_p\ot_{Z_p^G} U^G$ as algebras.
The first two of these isomorphisms are known as
Friedlander-Parshall factorisations; see \cite{FrPa} Theorem 4.1.

\smallskip
\noindent
2. It also follows from Proposition~2
that $Z$ is integral over $Z_p$ (for $Z$ is a finitely generated $Z_p$-module and
$Z_p$ is Noetherian\footnote{This also follows from
a version of  the PWB theorem; see
\cite{Jac}, Chap. 5, \S~7, Lemma 4.}). So
${\mc Q}(\g)$ is a
finite extension of the field of fractions of
$Z_p\cong S(\g)^{(1)}$ and hence $\mathrm{tr.\,deg}_K\,{\mc Q}(\g)=
\dim \g$. The analogous statements for the fields of fractions of
$K[\g]^\g$ and $S(\g)^\g$ are obvious.
\section{Proof of the main theorems}\label{thetheorems}
\newcommand{\e}{\varepsilon}
\noindent
{\bf 4.1.}\quad
Define
$\partial_{ij}\in\Der_K\,K[\gl_n]$ be setting
$\partial_{ij}(\xi_{rs})=1$ if $(r,s)=(i,j)$ and $0$ otherwise.
It is immediate from our discussion in 2.4 that
$s_k$ is the sum of all diagonal minors of order $k$ of the matrix
$\sum_{i,j}\xi_{ij}\,e_{i,j}$ with entries in $K[\gl_n]$.
If we write each $s_k$ as a polynomial in the $\xi_{ij}$, then we obtain $n$ equations in
the $\xi_{ij}$ and the $s_k$. By the above, $\xi_{ij}$ with one fixed row or column index
are not multiplied among each other in these equations. In particular these equations are
{\it linear} in $\xi_{11},\xi_{12},\ldots,\xi_{1n}$.

Let $R$ denote the ${\mb F}_p$-subalgebra of $K[\gl_n]$ generated by
all $\xi_{ij}$ with $i>1$.  Set $$M=
\left(\begin{array}{cccc}
\partial_{11} (s_1)&\partial_{12}(s_1)&\ldots& \partial_{1n}(s_1)\\
\partial_{11} (s_2)&\partial_{12}(s_2)&\ldots& \partial_{1n}(s_2)\\
\vdots&\vdots& &\vdots \\
\partial_{11} (s_n)&\partial_{12}(s_n)&\ldots& \partial_{1n}(s_n)
\end{array}\right),\quad
{\bf c}=\left(\begin{array}{c}
\xi_{11}\\
\xi_{12}\\
\vdots \\
\xi_{1n}
\end{array}\right),\quad
{\bf s}=\left(\begin{array}{c}
s_1\\
s_2\\
\vdots \\
s_n
\end{array}\right).
$$
By the preceding paragraph the matrix $M$ has entries in $R$ and the following
vector equation holds:

\begin{equation}
M\cdot {\bf c}\,=\,{\bf s} + {\bf r},  \quad\,\mbox{where} \  \
M\in\gl_n(R)\ \, \mbox{and}\ \ {\mathbf r}\in R^n.
\end{equation}
Clearly, $M$ is a matrix with functional entries.
Hence $M(x)\in\gl_n$ is well-defined
for any $x\in \gl_n$.
Let $d=\det M$, a regular function on $\g$. Recall from 3.1
the definition of the affine subspace
${\mc S}=\{x_{\mathbf a}\,|\,{\bf a}\in K^n\}$ of $\gl_n$.
\begin{lemgl}
For all ${\bf a}\in K^n$ we have
$d(x_{\bf a})=(-1)^{\lfloor n/2\rfloor}$. In particular, $d\ne 0$.
\end{lemgl}
\begin{proof}
Let $\xi_1,\ldots, \xi_n$ be the coordinate functions on $K^n$ and let $\partial_i$ be
the derivation of the coordinate ring of
$K^n$ such that $\partial_i(\xi_j)=1$ when $i=j$ and $0$ otherwise.
Then it is easy to see that $\partial_{1j}(f)(x_{\bf a})=
\partial_j\big({\bf b}\mapsto f(x_{\bf b})\big)({\bf a})$ for all $f\in K[\gl_n]$.
Furthermore it follows from the formula displayed in 2.4 and our remarks in the
proof of Lemma~1 that $s_i(x_{\bf a})=(-1)^{i-1}a_i$.
So the $(i,j)$th entry of $M(x_{\bf a})$ equals $(-1)^{i-1}\partial_j(\xi_i)$.
But then $M(x_{\bf a})=\text{diag}(1,-1, \ldots,(-1)^{n-1})$ and
$\det(M)(x_{\bf a})=(-1)^{\lfloor n/2\rfloor}$.
\end{proof}

\smallskip
\noindent
{\bf 4.2.}\quad Let $Q$ denote the field of fractions of
$K[\g]^{\g}$.
It follows from Proposition~\ref{prop.g-invariants} that
$Q$ is generated by $m+\dim \g$ elements.
Using Lemma~3 we will show that $m$ generators can be made
redundant here. Since
$\mathrm{tr.\,deg}_K\,Q=\dim \g$, this will imply that $Q$ is rational.
We will then use a very similar method
to establish the rationality of ${\mc Q}$.

Let $F\colon\,f\mapsto f^p$ denote the Frobenius endomorphism of $K[\gl_n]$.
It acts componentwise on $\gl_n\big(K[\gl_n]\big)$ and
$K[\gl_n]^n$. Note that $R^F\subset R$.
\begin{thmgl}\label{thm.rat}
Both $S(\g)^\g$
and $Z$ have rational fields of fractions.
\end{thmgl}
\begin{proof}
1. First we assume that $\g=\gl_n$. Applying $F$ to both sides of
(1) we get
\begin{equation}
M^F\cdot\, {\bf c}^F\,=\,{\bf s}^F + {\bf r}^F,  \quad\,\mbox{where} \  \,
M\in\gl_n(R^p)\ \ \mbox{and}\ \ {\mathbf r}\in (R^p)^n.
\end{equation}
By Lemma~3, $\det(M^F)=d^p\ne 0$. Therefore, ${\bf c}^F$ has components in
the ${\mb F}_p$-subalgebra of $Q$ generated by $
s_1^p,\ldots, s_n^p$, $\,(d^{p})^{-1}$ and $\xi_{ij}^p$ with $i>1$.
As a result, $Q$
is generated  by the $n^2$ elements
$s_1,\ldots, s_n$ and $\xi_{ij}^p$ with $i>1$. These elements must be algebraically
independent because $\mathrm{tr.\,deg}_K\,Q=n^2$;  see Remark~2.
Thus $Q$ is rational over $K$.
The same assertion then holds for the
field of
fractions of $S(\g)^\g$ in view of the $G$-equivariant algebra
isomorphism $\theta\colon\,K[\g]\,\stackrel{\sim}{\longrightarrow}\, S(\g)$.

\smallskip

\noindent
2. Recall from 2.2 and 2.4 that $\eta\circ\theta:\, K[\g]^{(1)}\to Z_p$
is a $G$-equivariant algebra isomorphism.
Observe that $\theta(\xi_{ij})=e_{j,i}$ and that
${\mc R}:=\eta(\theta(R))$ is the ${\mb F}_p$-subalgebra
of $Z_p$ generated by all $e_{i,j}^p-e_{i,j}^{[p]}$ with $j>1$. Let $\mathbf e\in Z_p^n$
denote
the column vector whose $i$th component equals $e_{i,1}^p-e_{i,1}^{[p]}$.
Applying $\eta\circ\theta$ to both sides of (2)  yields
\begin{equation}
{\mc M}\cdot {\bf e}=\eta(\theta({\bf s})) + \tilde{{\bf r}}, \quad\,\mbox{where} \  \,
{\mc M}\in\gl_n({\mc R})\ \ \mbox{and}\ \ \tilde{{\mathbf r}}\in {\mc R}^n.
\end{equation}
By
Proposition~\ref{prop.g-invariants},  ${\mc Q}$ is generated
over $K$ by the
elements $e_{i,j}^p-e_{i,j}^{[p]}$ and $n$ algebraically independent elements
generating $Z^G$. Besides, $\mathrm{tr.\,deg}_K\,{\mc Q}=n^2$;
see Remark~2.
Since
$\eta(\theta(s_i))\in Z_p^{G}$ and $\det {\mc M}=\eta(\theta(d))\ne 0$, we now  argue
as in part~1 of this proof  to deduce that ${\mc Q}$ is rational over $K$.

\smallskip

\noindent
3. Now assume that $\g=\spl_n$ and $p\nmid n$.
Recall our notation from 2.4.
Applying the restriction
homomorphism $K[\gl_n]\to K[\spl_n],\ f\mapsto f',$ to both sides of
(1) we obtain the following equations in the $\xi_{ij}'$ and $s_2',\ldots,s_n'$
$$M'\cdot \,{\bf c}'={\bf s}' + {\bf r}',\quad\,\mbox{where} \  \,
M'\in\gl_n(R')\ \ \mbox{and}\ \ {\mathbf r}'\in (R')^n.$$ Here $M'$, ${\bf c}'$,
${\bf s}'$, ${\bf r}'$
have the obvious meaning
and $R'$ is the $\mb{F}_p$-subalgebra of $K[\spl_n]$
generated by all $\xi'_{ij}$ with $i>1$.
Note that we now have
$\theta(\xi_{ij}')=e_{j,i}$ for $i\neq j$ and $\theta(\xi_{ii}')=e_{i,i}-(1/n)I_n$.
Since
${\mc S}\cap\spl_n\ne\emptyset$, Lemma~3 shows that
$d'=\det( M')\ne 0$. We can  thus repeat our
arguments from parts~1 and 2 of this proof
to deduce that the generators $(\xi_{11}')^p,\ldots, (\xi_{1n}')^p$ and
$\big(e_{1,1}-(1/n)I_n\big)^p-
\big(e_{1,1}-(1/n)I_n\big),\,e_{2,1}^p,\ldots, e_{n,1}^p$ of $Q$ and ${\mc Q}$,
respectively,
are redundant. This proves that both $Q$ and $\mc Q$ are
$K$-rational in the present case (recall that we now have one
generator
less  and $\mathrm{tr.\,deg}_K\,Q=\mathrm{tr.\,deg}_K\,{\mc Q}=
n^2-1$; see Proposition~\ref{prop.g-invariants} and Remark~2).  \end{proof}
\noindent
{\bf 4.3.}\quad
We now turn our attention to the second problem: the unique factorisation property.
The determinant $d$ will play a prominent r{\^o}le here.
\begin{propgl}
The polynomial function $d$ is irreducible in $K[\gl_n]$.
\end{propgl}
\begin{proof}
1. Let $\g=\gl_n$ and
let $P$ be the maximal parabolic subgroup of $G=GL_n(K)$
consisting of all invertible
matrices $(\lambda_{ij})$ with $\lambda_{i1}=0$ for all $i>1$. As a first step,
we are going to show that
$d$ is a semiinvariant for $P$. We have
\begin{equation}
d\,=\sum_{\pi\in {\mf S}_n}
\sgn(\pi)\,\partial_{1,\pi(1)}(s_1)\,\cdots\,\partial_{1,\pi(n)}(s_n).
\end{equation}
The adjoint action of $G$ on $\g$ induces a natural action of $G$
on the Lie algebra $\Der_K\,K[\g]$.
The subspace $\mf D$ of $\Der_K\,K[\g]$
consisting of all homogeneous  derivations of degree
$-1$  is $G$-stable and has $\{\partial_{ij}\,|\,\,1\le i,j\le n\}$
as a basis.
We define ${\mf D}_0$ to be the subspace of $\mf D$ spanned by all
$\partial_{1i}$ with $1\le i\le n$.

Let $\g_0^*$ denote the subspace $\g^*$ spanned by all  $\xi_{i,j}$ with
$i>1$. It is easy to see
that $\g_0^*$ consists of all linear functions $\psi$ on $\g^*$ with
$\psi(e_{1,i})=0$ for all $i$.
As the linear span of  all $e_{1,i}$  is
$(\Ad\,P)$-invariant, $\g_0^*$ is invariant under the
coadjoint action of $P$ on $\g^*$. As ${\mf D}_0=\{D\in{\mf D}\,|\,\,
\g_0^*\subset \Ker\,D\}$, it follows that
$g\circ D\circ g^{-1}\in{\mf D}_0$ for all $D\in{\mf D}_0$ and
$g\in P$. Thus $P$ acts on ${\mf D}_0$. We denote by
$\tau$ the corresponding representation of $P$.

Let $g$ be any element in $P$ and denote by $A=(a_{ij})$  the
matrix of $\tau(g)$  relative to the
basis
$\{\partial_{1i}\,|\,\,1\le i\le n\}$ of ${\mf D}_0$.
Since each $s_i$ is $G$-invariant, we have
$$g\big(\partial_{1j}(s_i)\big)=
(g\circ\partial_{1j}\circ g^{-1})(s_i)=\big(\tau(g)(\partial_{1j})\big)(s_i)
\qquad\,(1\le i,j\le n).$$
Combining this with (4) we now obtain
\begin{eqnarray*}
g(d)&=&\sum_{\pi\in {\mf S}_n}
\sgn(\pi)\,\big(\tau(g)(\partial_{1,\pi(1)})\big)(s_1)\,\cdots\,
\big(\tau(g)(\partial_{1,\pi(n)})\big)(s_n)\\
&=&
\sum_{\pi\in {\mf S}_n}
\sgn(\pi)\,\big(\textstyle{\sum}_k\,a_{k,\pi(1)}\,\partial_{1,k}(s_1)\big)\,\cdots\,
\big(\textstyle{\sum}_k\,a_{k,\pi(n)}\,\partial_{1,k}(s_n)\big)\\
&=&\det\big(\big(\textstyle{\sum}_k a_{kj}\,\partial_{1k}(s_i)\big)_{ij}\big)
=\, \det(M\cdot A)\,=\,(\det A)\,d.
\end{eqnarray*}

\smallskip

\noindent
2. Let $B$ be the Borel subgroup of $G$ consisting of all invertible upper
triangular matrices. Clearly, $B\subset P$. Since $d$ is a semiinvariant for $P$,
the Borel subgroup $B$ acts on the line $K d$ through a rational character,
say $\chi$.
Let $T$ be as in 2.3, a maximal torus of $G$ contained in $B$.
We need to determine the weight of $d$ with respect to $T$.
Note that the maximal unipotent subgroup $U^+$ of $B$ acts trivially on $K d$.

Let $X(T)$ denote the lattice of rational characters of $T$.
For $i\in\{1,\ldots,n\}$
we denote by $\e_i$ the rational character $\text{diag}\,(\lambda_1,\ldots,
\lambda_n)\mapsto \lambda_i$ of $T$.
It is well-known that $X(T)$ is a free $\mb Z$-module with $\e_1,\ldots,\e_n$ as a
basis, and
$\Sigma=\{\e_i-\e_j\,|\,i\neq j\}$ is the set of roots of $G$ with respect to $T$.
For $1\le i\le n-1$ put $\alpha_i=\e_i-\e_{i+1}$.
It is well-known that $\Sigma$
is a root system of type $A_{n-1}$ in its $\mb{R}$-span in
$\mb{R}\ot_\mb{Z}X(T)$ and, moreover,
$\alpha_1,\ldots,\alpha_{n-1}$ form
the basis of simple roots of
$\Sigma$ relative to $B$. We denote the
corresponding fundamental weights
by $\varpi_1,\ldots,\varpi_{n-1}$.

From the fact that $\xi_{ij}$
has weight $\e_j-\e_i$ relative to $T$  it follows
that $\partial_{ij}$ has weight $\e_i-\e_j$. This implies that
all summands
in (4) have the same $T$-weight
$\sum_{i=1}^n(\e_1-\e_i)=n\e_1-\sum_{i=1}^n\e_i$ which is therefore also the $T$-weight
of $d$. Using Bourbaki's tables
it is now easy to observe that $\chi|_T=n\varpi_1$; see \cite{Bou1}.

\smallskip

\noindent
3. Now we will show that $d$ is irreducible. Let $d=f_1^{m_1}\cdots\, f_r^{m_r}$
be the prime factorisation of $d$ in the factorial ring
$K[\g]$. Since $d$ is homogeneous, so are
all $f_i$. By the
uniqueness of prime factorisation, the group $B$
permutes the lines $Kf_1,\ldots,Kf_r$.
Since $B$ is connected,
each $f_i$ is a semiinvariant for $B$.
Let $\chi_i$ denote the character of $B$ through which $B$
acts on $Kf_i$.

Observe that all weights of the $G$-module
$K[\g]$ are in the root lattice of $\Sigma$.
Since $U^+$ fixes $f_i$,  it must be that
$\chi_i|_T=\sum_{j=1}^{n-1}k_{i,j}\,\varpi_j$ where all $k_{i,j}$
are nonnegative integers;
see e.g. Proposition~II.2.6 in \cite{Jan1}.
The prime factorisation of $d$ and the concluding remark in part~2 of this proof
yield $$n\varpi_1\,=\,\,\sum_{i=1}^rm_i\big(\sum_{j=1}^{n-1}k_{i,j}\,
\varpi_j\big)\,=\,\,\sum_{j=1}^{n-1}\big(\sum_{i=1}^rm_i k_{i,j})\,\varpi_j.$$
Since all $m_i$ are strictly positive, we obtain that $n=\sum_{i=1}^rm_ik_{i,1}$ and
$k_{i,j}=0$ for all $j>1$. Since all
$\chi_i|_T=k_{i,1}\,\varpi_1$
are in the root lattice of $\Sigma$,
it must be that $n\!\mid\! k_{i,1}$ for all $i$. So there is $j$ such that
$k_{j,1}=n,$ $m_j=1$
and  $k_{i,1}=0$ for $i\ne j$.
In other words,
$d=d_1 d_2$ where $d_1$ is an {\it irreducible}
semiinvariant for $B$ and $d_2$ is a {\it homogeneous}
polynomial function on $\g$  invariant under $T\,U^+=B$.

On the other hand, it is well-known that $K[\g]^B=\,K[\g]^G$ (this is immediate
from the completeness of the  flag variety $G/B$).
Hence $d_2\in K[s_1,\ldots,s_n]$. Since $s_i(x_{\bf 0})=0$ for all $i$,
Lemma~3 shows that $d_2$ is a nonzero scalar.
We conclude that $d$ is irreducible as desired.
\end{proof}
\begin{cornn}
If $p\nmid n$, then the polynomial function $d'$ is irreducible in $K[\spl_n]$.
\end{cornn}
\begin{proof}
Let $G=\GL_n(K)$. The restriction map
$K[\gl_n]\to K[\spl_n]$  is $G$-equivariant.
As in  parts~1 and 2 of the previous proof one proves that
$d'$ is a semiinvariant for $P$ of weight $n\varpi_1$. The argument in part~3
then shows that $d'$ is irreducible.
\end{proof}
\noindent
{\bf 4.4.} We will need a result from Commutative Algebra
often referred to as {\it Nagata's lemma}; see  \cite{Eis}, Lemma 19.20, for
example. It asserts the following:
If $x$ is a prime element of a Noetherian integral domain $S$
such that  $S[x^{-1}]$ is factorial,  then $S$ is factorial.
\begin{thmgl}\label{thm.UFD}
The centre of $U(\g)$ is a  unique factorisation domain.
\end{thmgl}
\begin{proof}
1. Suppose $\g=\gl_n$, where $n\ge 2$, and set
$d_0=\eta(\theta(d))$. It is immediate from (3) that  $Z[d_0^{-1}]$ is
isomorphic to a localisation of a polynomial algebra in $\dim \g$ variables.
Since any localisation of a factorial ring is again factorial,
$Z[d_0^{-1}]$ is is a unique factorisation domain.
We claim that $d_0$ is a prime element of $Z$.
Our remarks in 2.1, 2.2 and 2.4 show that $\gr\,(d_0)=\theta(d^p)$ and that
$$\gr\big(Z/(d_0)\big)\cong S(\g)^\g/\big(\theta(d^p)\big)\cong K[\g]^\g/(d^p).$$
Hence the claim will follow if we establish that
$K[\g]^\g/(d^p)$ has no zero divisors; see 2.1 for more detail.

By our remarks in  the proof of Proposition~3,
the semiinvariant
$d$ has weight $\chi|_T=n\varpi_1$ relative to $T$. So
$\chi|_T\not\in p\,X(T)$, for $n\varpi_1$ is indivisible in $X(T)$.
It follows that
the Lie algebra $\h={\Lie}\,T$ does not
annihilate $d$. As a result, $d\not\in K[\g]^\g$.
So Proposition~3 and Lemma~2 yield
that $d^p$ is an irreducible element
of the factorial ring $K[\g]^\g$.
But then $K[\g]^\g/(d^p)$ has no zero divisors, as wanted.
Thus $d_0$ is a prime element of $Z$. Applying Nagata's lemma we finally deduce
that $Z$ is factorial in the present case.

\smallskip

\noindent
2. Suppose $\g=\spl_n$ and $p\nmid n$. If $(p,n)=(2,2)$, the factoriality of $Z$
has  been established in \cite{BHa}, so assume that this is not the case.
Then the Lie algebra $\g$ is perfect. Put $T_0=T\cap\SL_n(K)$. The restriction homomorphism
$X(T)\to X(T_0)$ induces an isomorphism of root systems. We denote the images of the $\alpha_i$
and $\varpi_i$ under this isomorphism by the same symbols.
Now the weight lattice of the root system $\Sigma$
coincides with the character group $X(T_0)$.
By the proof of Proposition~3 $d'$ has weight
$n\varpi_1$ relative to $T_0$. Since $p\nmid n$,
we have  $n\varpi_1\not\in p\,X(T_0)$.
So the Lie algebra $\h'={\Lie}\,T_0$
does not annihilate $d'$, forcing $d'\not\in K[\g]^\g$.
In view of Corollary~2 and Proposition~1 this shows that
$(d')^p$ is a prime
element of the factorial ring $K[\g]^\g$.

Now set $d'_0=\eta(\theta(d'))$.
Repeating the argument from the beginning of part~1 of this
proof we now see that $d'_0$ is a prime element of $Z$. A version of (3) for
$\g=\spl_n$ with $p\nmid n$
implies that $Z[(d_0')^{-1}]$ is a unique factorisation domain. But then so is
$Z$, by Nagata's lemma, completing the proof.

\end{proof}

\end{document}